\documentclass[12pt]{article}
\usepackage{}
\usepackage{mathrsfs,graphicx,amsmath,amssymb,amsfonts}

\usepackage{appendix}
\usepackage{caption}
\usepackage{mathrsfs,graphicx,amsmath,amssymb,amsfonts,amsthm}
\usepackage{cite}
\usepackage{longtable}
\usepackage{supertabular}
\evensidemargin=\oddsidemargin\topmargin=-1.5cm

\theoremstyle{definition}

\begin{document}
\title{The least signless Laplacian eigenvalue of the complements of bicyclic graphs \thanks{This work is supported by NSFC (No. 11461071).}}
\author{Xiaoyun Feng, Guoping Wang\footnote{Corresponding author. Email: xj.wgp@163.com.}\\
{\small School of Mathematical Sciences, Xinjiang Normal University,}\\
{\small Urumqi, Xinjiang 830054, P.R.China}}
\date{}
\maketitle {\bf Abstract.}
Suppose that $G$ is a connected simple graph with the vertex set $V(G)=\{v_1, v_2,\cdots,v_n\}$.
Then the adjacency matrix of $G$ is $A(G)=(a_{ij})_{n\times n}$,
where $a_{ij}=1$ if $v_i$ is adjacent to $v_j$, and otherwise $a_{ij}=0$.
The degree matrix $D(G)=diag(d_{G}(v_1), d_{G}(v_2), \dots, d_{G}(v_n)),$
where $d_{G}(v_i)$ denotes the degree of $v_i$ in the graph $G$ ($1\leq i\leq n$).
The matrix $Q(G)=D(G)+A(G)$ is called the signless Laplacian matrix of $G$.
The least eigenvalue of $Q(G)$ is also called the least signless Laplacian eigenvalue of $G$.
In this paper we give two graft transformations and then use them to characterize the unique connected graph
whose least signless Laplacian eigenvalue is minimum among the complements of all bicyclic graphs.

{\flushleft{\bf CLC number:}}  O 157.5

\section{Introduction}
~~~~Suppose that $G$ is a connected simple graph with the vertex set $V(G)=\{v_1, v_2,\cdots,v_n\}$.
Then the adjacency matrix of $G$ is $A(G)=(a_{ij})_{n\times n}$,
where $a_{ij}=1$ if $v_i$ is adjacent to $v_j$, and otherwise $a_{ij}=0$.
The degree matrix $D(G)=diag(d_{G}(v_1), d_{G}(v_2), \dots, d_{G}(v_n)),$
where $d_{G}(v_i)$ denotes the degree of $v_i$ in the graph $G$ ($1\leq i\leq n$).
The matrix $Q(G)=D(G)+A(G)$ is called the signless Laplacian matrix of $G$.
Since $Q(G)$ is positive semidefinite,
its eigenvalues can be arranged as $\lambda_1(G)\geq \lambda_2(G)\geq \dots \geq \lambda_n(G)\geq 0$,
where $\lambda_n(G)$ is also called the least signless Laplacian eigenvalue of $G$, denoted by $\lambda(G)$.

The least signless Laplacian eigenvalues of connected graphs have been studied extensively.
Cardoso, Cvetkovi\'c, Rowlinson and Simi\'c \cite{Cardoso} determined the unique graph
whose least signless Laplacian eigenvalue attains the minimum among all connected non-bipartite graphs.
Guo, Chen and Yu \cite{SGGuo} obtained a lower bound on the least signless Laplacian eigenvalue of a graph.
Guo and Zhang \cite{SGGGuo} described the unique graph whose least signless Laplacian eigenvalue attains the minimum
among all non-bipartite connected graphs with fixed maximum degree.
Fan, Wang and Guo \cite{Fan} determined the graph whose least signless Laplacian eigenvalue
attains the minimum or maximum among all connected non-bipartite graphs with fixed order and given number of pendant vertices.
Guo, Ren and Shi\cite{JMGuo} determined the graph whose the least signless Laplacian eigenvalue
attains the maximum among all connected unicyclic graphs.
He and Zhou \cite{CXHe} gave a sharp upper bound on the least signless Laplacian eigenvalue of a graph using domination number.
Wang and Fan \cite{YWang} determined the graph whose the least signless Laplacian eigenvalue is minimum.
Yu, Fan and Wang \cite{GDY} determined the unique graph
whose least signless Laplacian eigenvalue attains the minimum among all connected non-bipartite bicyclic graphs.
Yu, Guo and Xu \cite{GYu} determined the unique graph whose least signless Laplacian eigenvalue attains the minimum
among all connected non-bipartite graphs with given matching number and edge cover number, respectively.
Wen, Zhao and Liu \cite{QWen} determined the graph which has the minimum the least signless Laplacian eigenvalue
among all non-bipartite graphs with given stability number and covering number, respectively.

Suppose that $G=(V(G),E(G))$ is a connected simple graph.
The complement of $G$ is denoted by $G^c=(V(G^c), E(G^c))$,
where $V(G^c)=V(G)$ and $E(G^c)=\{xy:x,y\in V(G), xy\notin E(G)\}.$
If $\mid E(G) \mid=\mid V(G) \mid+1$ then $G$ is {\it bicyclic graph}.
Denote by $K_{1,n-1}$ the star graph on n vertices,
and by $k_{1,n-1}+2e$ the graph which is obtained from $k_{1,n-1}$ by connecting two pairs of different pendant vertices.
The complement $(k_{1,n-1}+2e)^c$ of $k_{1,n-1}+2e$ contains an isolated vertex,
and so it is not a connected graph.
We also note that the complement of arbitrary bicyclic graph is non-bipartite graph if its order is greater or equal $12$.
Therefore, we will only consider the complements of those connected bicyclic graphs on $n\geq 12$ vertices except $k_{1,n-1}+2e$.

Li and Wang \cite{SCLi} determined the unique graph whose least signless Laplacian eigenvalue attains the minimum
in the set of the complements of all trees except $K_{1,n-1}$.
Yu, Fan and Ye \cite{GDYu} obtained the unique graph whose least signless Laplacian eigenvalue
attains the minimum in the set of the complements of all unicyclic graphs except $K_{1,n-1}+e$.
In this paper we will give two graft transformations and then use them to characterize the unique graph
whose least signless Laplacian eigenvalue is minimum
among the complements of all bicyclic graphs on $n\geq 12$ vertices except $k_{1,n-1}+2e$.

\section {Main results}
Suppose that $G$ is a graph with the vertex set $V(G)=\{v_1, v_2,\ldots,v_n\}$.
Let $X=(X_1, X_2, \ldots, X_n)^T$ be an unit vector such that $X_i=X(v_i)$ ($1\leq i\leq n$).
Then we have
$$ X^{T}Q(G)X= \sum\limits_{{v_i}{v_j}\in E(G)}(X_i+X_j )^2 ~~\eqno{(2.1)}$$

$$\lambda(G)\leq X^{T}Q(G)X  \eqno {(2.2)}$$
The equality (2.2) holds if and only if $X$ is the eigenvector of $Q(G)$ corresponding to $\lambda(G)$.\vskip 2mm

In what follows the vertex sets of all graphs on $n$ vertices are write as
$\{v_1, v_2,\ldots,v_n\}$, and let $X=(X_1, X_2, \ldots, X_n)^T$ be an unit vector such that $X_i=X(v_i)$ ($1\leq i\leq n$)
and $X_1 \geq X_2 \geq\dots\geq 0\geq\dots\geq X_n$. \vskip 2mm

{\bf Lemma 2.1.} \cite{SCLi}. {\it Let $X$ be as above with $X_1> 0$ and $X_n< 0$.
Then for any $1\leq i,j \leq n$,
$(X_{i}+ X_{j} )^2 \leq max  \{(X_{i}+ X_{1} )^2,(X_{i}+ X_{n} )^2 \}$ and $(X_{i}+ X_{j} )^2 \leq max\{(X_{j}+ X_{1} )^2,(X_{j}+ X_{n} )^2 \}.$} \vskip 2mm

Let $\mathcal{C}_n^c$  $(n\geq 12)$ be the set of all connected graphs
each of which is the complement of a connected bicyclic graph on $n$ vertices.
Then, for any $G \in \mathcal{C}_n^c$, $\lambda(G)>0$. \vskip 2mm

Suppose that $v_i$ and $v_j$ are two vertices of a graph $G$.
Then the {\it distance} $d_G(v_i,v_j)$ between $v_i$ and $v_j$ is the length of the shortest path between $v_i$ and $v_j$. \vskip 2mm

Define a {\it $b$-graph} to be a graph consisting of two vertex-disjoint
cycles and a path joining them having only its
end-vertices in common with the cycles. Define a
{\it $\infty$-graph} to be a graph consisting of two cycles
with exactly one vertex in common. Define a
{\it $\theta$-graph} to be a graph consisting of two basic cycles with at least two vertices in common.
Obviously, a bicyclic graph is one of $b$-graph, $\infty$-graph and $\theta$-graph with trees attached.

Let the graphs $G_k{(p,q)}$ ($1\leq k\leq 12$) be shown in Fig. 1.
\begin{center}
\begin{figure}[htbp]
\centering
\includegraphics[height=45mm,  width=130mm]{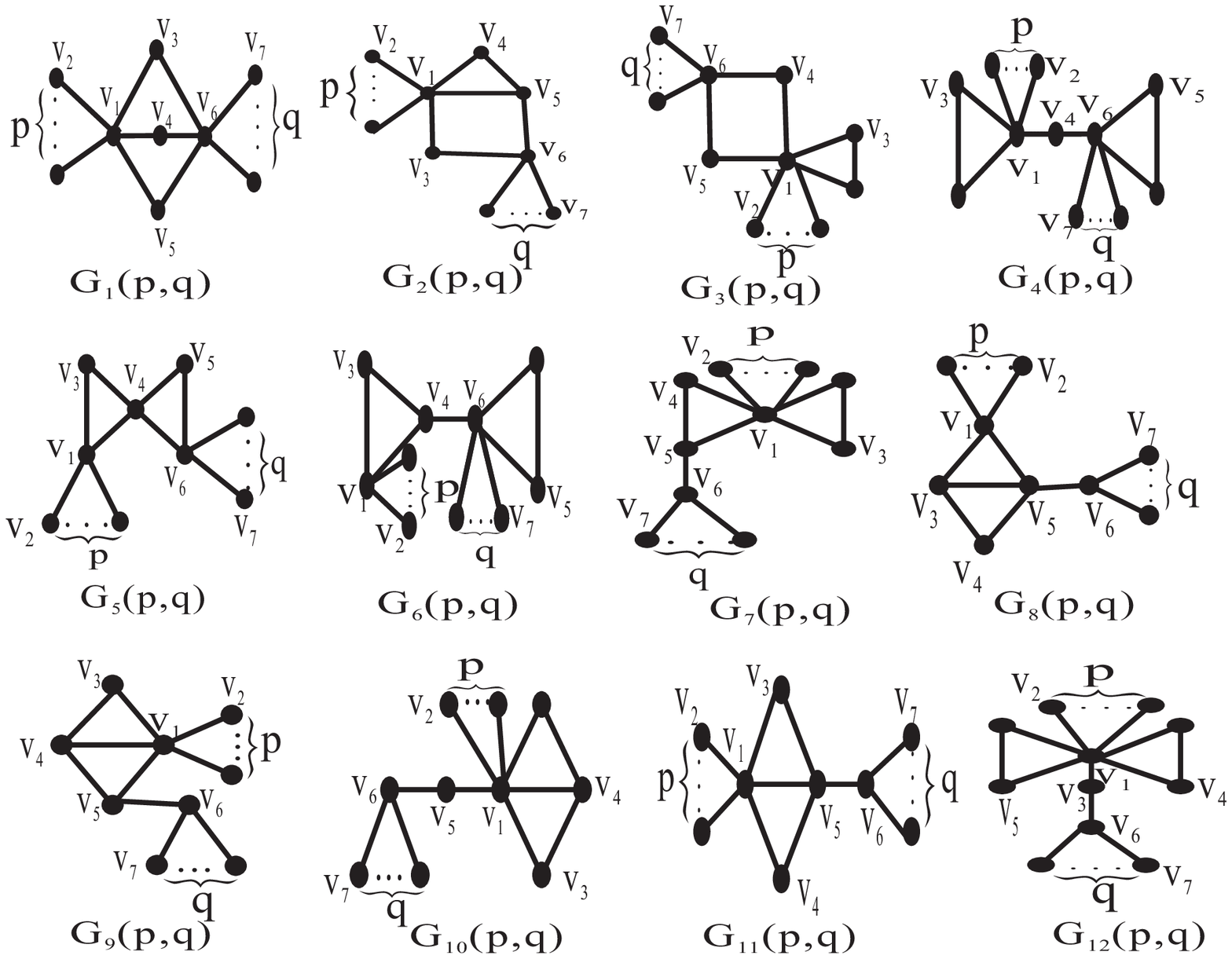}
\\{Fig. 1. $G_{k}(p,q)$ ($1\leq k\leq 12$)}
\end{figure}
\end{center}

Suppose that $G$ is a graph on $n$ vertices and
that $v_1v_{l_1}v_{l_2}\cdots v_{l_t}v_n$ is the shortest path between $v_1$ and $v_n$.
Let $X=(X_1, X_2, \ldots, X_n)^T$ satisfy $X_1> 0$ and $X_n< 0$.
Delete $v_{l_1}v_{l_2}$, and add the edge $v_1v_{l_2}$ if $(X_1+X_{l_2})^2\geq (X_n+X_{l_1})^2$,
and otherwise add the edge $v_nv_{l_1}$.
The above process is called {\it $T_1$-transformation} of $G$.
Next we always denote by $G_{T_1}$ the result graph which is obtained after some steps of $T_1$-transformation for $G$.
Now we use $T_1$-transformation to prove the below result. \vskip 2mm

{\bf Lemma 2.2.} {\it Let $\mathcal{H}^{'}=\{G_k{(p,q)}\mid 1\leq k\leq 12\}$.
If $G^c\in \mathcal{C}_n^c$ then either $G_{T_1}\in \mathcal{H}^{'}$ or $d_{G_{T_1}}(v_1,v_n) \leq 2$.}

{\bf Proof.} If $G \in \mathcal{H}^{'}$ or $d_{G}(v_1,v_ n)\leq 2$ then  it is finished.
So we next assume $G\notin \mathcal{H}^{'}$ and $d_{G}(v_1,v_n) >2$.
Now we distinguish four cases to discuss.

{\bf Case 1.} $v_1$ and $v_n$ are on the same cycle of $G$.

Making some steps of $T_1$-transformation for $G$ we can see that $d_{G_{T_1}}(v_1,v_n)\leq d_{G}(v_1,v_n)-1$.
Hence we can determine that the case is true.

{\bf Case 2.} $v_1$ and $v_n$ are on the different cycles of $G$.

In this case $G$ is the $b$-graph or $\infty$-graph with trees attached.
Making some steps of $T_1$-transformation for $G$
we can observe that either $G_{T_1}$ is $\theta$-graph with trees attached
or $d_{G_{T_1}}(v_1,v_n)\leq d_{G}(v_1,v_n)-1$.
If the former rises then the Case 1 shows that the case is true,
and otherwise repeat this process.
Thus we can determine that the case is true.

{\bf Case 3.} One of $v_1$ and $v_n$ is on the cycle and another is on the tree of $G$.

Suppose without loss of generality that $v_1$ is on the tree and that $v_n$ is on the cycle.
After making some steps of $T_1$-transformation for $G$
we can see that either $v_1$ and $v_n$ is on the same cycle of $G_{T_1}$
or $d_{G_{T_1}}(v_1,v_n)\leq d_{G}(v_1,v_n)-1$.
If the former rises then the Case 1 shows that the case is true,
and otherwise repeat this process.
Thus we can determine that the case is true.

{\bf Case 4.}  Both $v_1$ and $v_n$ are on the tree of $G$.

After making some steps of $T_1$-transformation for $G$
we can see that either one of $v_1$ and $v_n$ is on the cycle and another is on the tree of $G_{T_1}$,
or $d_{G_{T_1}}(v_1,v_n)\leq d_{G}(v_1,v_n)-1$.
If the former rises then the Case 3 shows that the case is true,
and otherwise repeat this process.
Hence we can determine that the case is true.  $\Box$ \vskip 2mm

If $X=(X_1, X_2, \ldots, X_n)^T$ satisfies $X^TQ(G)X=\lambda(G)$
then $X$ is the {\it unit first signless Laplacian eigenvector} of $G$.
Clearly,  for each $1\leq i\leq n$, we have
$$(\lambda(G)-d_{G}(v_i))X_i= \sum\limits_{v_j\in N_{G}(v_i)}X_j  \eqno{(2.3)}$$
where $N_{G}(v_i)$ denotes the neighbour set of $v_i$ in the graph $G$.
The equation (2.3) is also called the signless Laplacian eigenvalue-equation of $G$.

If $X=(X_1, X_2, \ldots, X_n)^T$ is an unit first signless Laplacian eigenvector of a graph $G$
such that $X_1 \geq X_2 \geq\dots\geq X_n$,
then $X_1> 0$ and $X_n< 0$
since the matrix $Q(G)-\lambda(G)I_n$ is positive semidefinite,
where $I_n$ is the identity matrix of order $n$.\vskip 2mm

We use $\theta_4$ and $\theta_5$ to denote respectively $\theta$-graph on $4$ and $5$ vertices,
and $\infty_5$ and $\infty_6$ to denote respectively $\infty$-graph on $5$ and $6$ vertices,
and $b_6$ and $b_7$ to denote respectively $b$-graph on $6$ and $7$ vertices whose two cycles are all length three.

Let the graphs $H_i$ $(1\leq i\leq 7)$ be shown in Fig. 2.
\begin{center}
\begin{figure}[htbp]
\centering
\includegraphics[height=33mm, width=110mm]{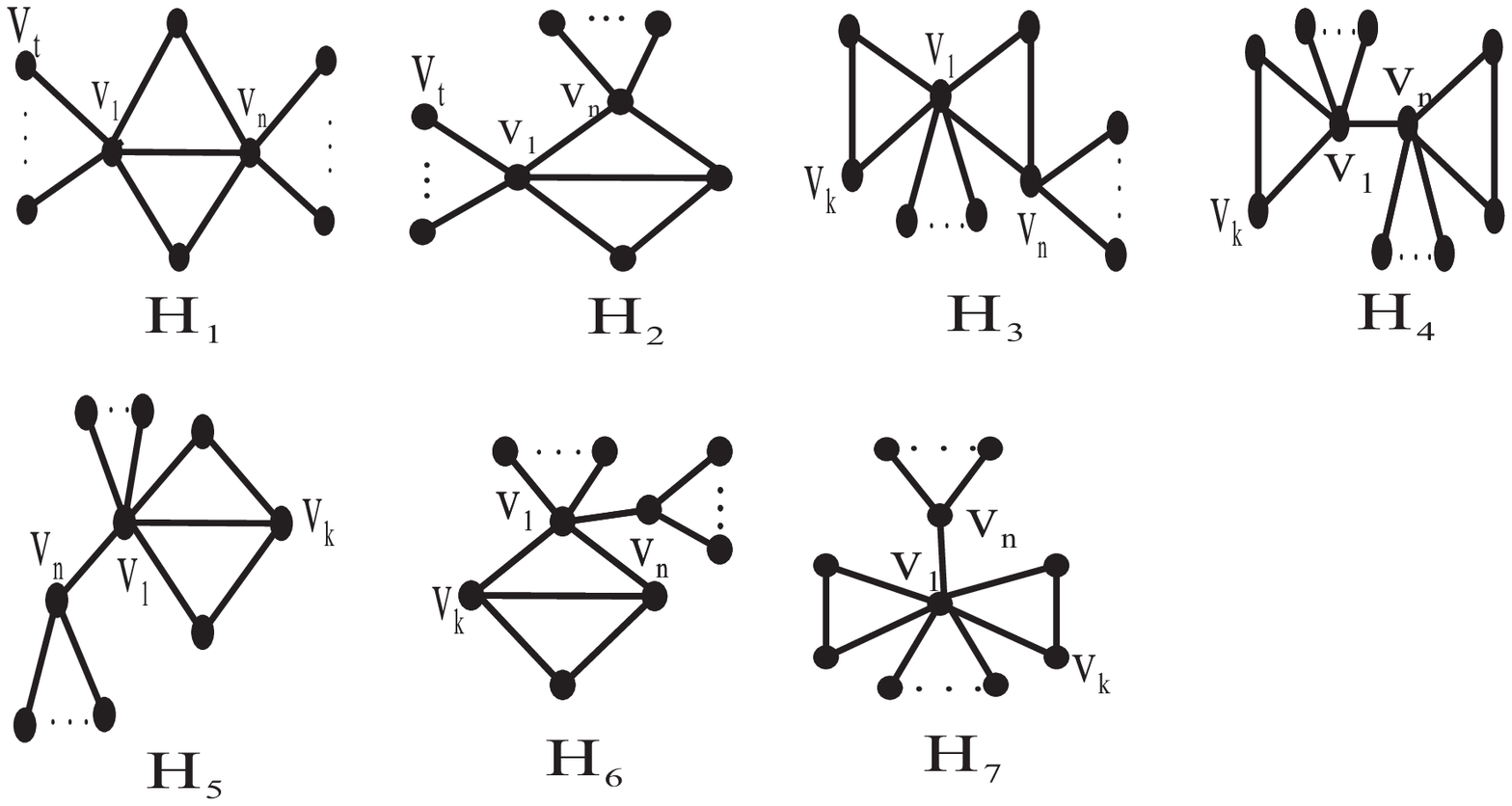}
\centering
\\{Fig. 2. $H_i$ $(1\leq i\leq 7)$}
\end{figure}
\end{center}

Suppose $H$ is the graph on $n$ vertices which is obtained from $\theta_4$, $\theta_5$, $\infty_5$,
$\infty_6$, $b_6$ or $b_7$ by attaching some trees to it.
Then we can find a pendent vertex $v_s$ whose neighbour $v_t$ is neither $v_1$ nor $v_n$.
Let $X=(X_1, X_2, \ldots, X_n)^T$ satisfy $X_1> 0$ and $X_n< 0$.
Delete $v_tv_s$, and add the edge $v_1v_s$ if $(X_1+X_{s} )^2\geq (X_n+X_{s})^2$,
and otherwise add the edge $v_nv_s$.
The above process is called {\it$T_2$-transformation} of $H$.
Next we always denote by $H_{T_2}$ the result graph which is obtained after some steps of $T_2$-transformation for $H$.
Now we use $T_1$ and $T_2$-transformations to prove the below important result. \vskip 2mm

{\bf Lemma 2.3.} {\it Let $\mathcal{H}^{''}=\{H_i\mid 1\leq i\leq 7 \}$.
Given a graph $G^c\in \mathcal {C}_n^c$ ($n\geq 12$),
there is a graph $H\in \mathcal{H}^{'}\cup \mathcal{H}^{''}$
such that $\lambda(G^c)\geq \lambda(H^c).$}

{\bf Proof.} Set $\mathcal{H}=\mathcal{H}^{'}\cup \mathcal{H}^{''}$
and $X=(X_1, X_2, \ldots, X_n)^T$ to be the unit first signless Laplacian eigenvector of $G^c$
such that $X_1 \geq X_2 \geq\dots\geq X_n$.
Then $X_1> 0$ and $X_n< 0$.
If $G \in \mathcal{H}$ then it is finished.
So we assume $G \notin \mathcal {H}$.

If $d_{G}(v_1,v_n) >2$
then making some steps of $T_1$-transformation for $G$,
by Lemma 2.2 we know that either $G_{T_1}\in \mathcal{H}$
or $G_{T_1}\notin \mathcal{H}$ but $d_{G_{T_1}}(v_1,v_n) \leq 2$.
If the former rises, then by the equation (2.1) and Lemma 2.1, we have
$$X^{T}Q(G)X= \sum\limits_{v_iv_j\in E(G)}(X_i+ X_j )^2\leq  \sum\limits_{v_iv_j\in E(G_{T_1})}(X_i+ X_j )^2=X^{T}Q(G_{T_1})X.$$
Now we assume the latter rises and distinguish two cases to discuss $G_{T_1}$.\vskip 2mm

\textbf{Case 1.} $d_{G_{T_1}}(v_1,v_n) = 1$.\vskip 2mm

\textbf{Case 1.1.} $v_1$ and $v_n$ are on the same cycle $C_{t+2}$ of order $t+2$.

If $t\geq 2$ then set $C_{t+2}=v_1v_n v_{l_1}v_{l_2}\cdots v_{l_t}v_1$.
Delete $v_{l_1}v_{l_2}$, and add the edge $v_1v_{l_1}$ if $(X_1+ X_{l_1} )^2\geq (X_n+ X_{l_2} )^2$, and otherwise add the edge $v_nv_{l_2}$.
Repeat this process until the result graph $\overline{G}_{11}$ contains the cycle $C_3 =v_1v_n v_{l_1} v_1$.
Set $v_{s_1}v_{s_2}\cdots v_{s_k}v_{s_1}$ to be another cycle which is contained in $\overline{G}_{11}$.
Delete $v_{s_1}v_{s_2}$, and add the edge $v_1v_{s_2}$ if $(X_1+X_{s_2} )^2\geq (X_n+X_{s_1} )^2$, and otherwise add the edge $v_nv_{s_1}$.
Repeating this process we can determine the result graph $\widetilde{G}$ is $\theta_4$ or $\infty_5$ with trees attached.
If $\widetilde{G}\notin \mathcal{H}$,
then making some steps of $T_2$-transformation for $\widetilde{G}$
we can get the graph $\widetilde{G}_{T_2}$ which is isomorphic to one of the graphs $H_1, H_2$ and $H_3$ as in Fig. 2.

\textbf{Case 1.2.}  $v_1$ and $v_n$ are on the different cycle.

In this case $G_{T_1}$ is the $b$-graph with trees attached.
Suppose $G_{T_1}$ contains the cycle $C_{t+1}=v_1v_{l_1}v_{l_2}\cdots v_{l_t}v_1$.
When $t\geq 3$ we delete $v_{l_1}v_{l_2}$, and add the edge $v_1v_{l_2}$ if $(X_1+X_{l_2} )^2\geq (X_n+X_{l_1} )^2$,
and otherwise add the edge $v_nv_{l_1}$.
Denote by $\overline{G}_{21}$ the result graph.
If $v_1$ and $v_n$ are on the same cycle of $\overline{G}_{21}$ then we enter Case 1.1,
and otherwise repeat this process until the result graph $\overline{G}_{22}$  contains the cycle $C_3 =v_1v_{l_1} v_{l_2} v_1$.
Next set $v_nv_{s_1}v_{s_2}\cdots v_{s_k}v_n$ to be another cycle which is contained in $\overline{G}_{22}$.
Delete $v_{s_1}v_{s_2}$,
and add the edge $v_nv_{s_2}$ if $(X_n+X_{s_2} )^2\geq (X_1+X_{s_1} )^2$, and otherwise add the edge $v_1v_{s_1}$.
Denote by $\overline{G}_{23}$ the result graph.
If $v_1$ and $v_n$ are on the same cycle of $\overline{G}_{23}$, then we enter Case 1.1,
and otherwise repeat this process until the result graph $\widetilde{G}$ is $b_6$ with trees attached.
If $\widetilde{G}\notin \mathcal{H}$,
then making some steps of $T_2$-transformation for $\widetilde{G}$.
We can determine that the graph $\widetilde{G}_{T_2}$ is isomorphic to the graph $H_4$ as in Fig. 2.

\textbf{Case 1.3.} One of $v_1$ and $v_n$ is on the cycle and another is on the tree.

Suppose without loss of generality that $v_1$ is on the cycle and that $v_n$ is on the tree
for otherwise we replace $X$ with $-X$.
Then $G_{T_1}$ contains the cycle $C_{t+1}=v_1v_{l_1}v_{l_2}\cdots v_{l_t}v_1$.
When $t\geq 3$ we delete $v_{l_1}v_{l_2}$, and add the edge $v_1v_{l_2}$ if $(X_1+X_{l_2} )^2\geq (X_n+X_{l_1} )^2$, and otherwise add the edge $v_nv_{l_1}$.
We denote by $\overline{G}_{31}$ the result graph.
If $v_1$ and $v_n$ are on the same cycle of $\overline{G}_{31}$ then we enter Case 1.1,
and otherwise repeat this process until the result graph $\overline{G}_{32}$ contains the cycle $C_3 =v_1v_{l_1} v_{l_2} v_1$.
Set $v_{s_1}v_{s_2}\cdots v_{s_k}v_{s_1}$ to be another cycle of $\overline{G}_{32}$.
Delete $v_{s_1}v_{s_2}$, and add the edge $v_1v_{s_2}$ if $(X_1+X_{s_2} )^2\geq (X_n+X_{s_1} )^2$, and otherwise add the edge $v_nv_{s_1}$.
We denote by $\overline{G}_{33}$ the result graph.
If $v_1$ and $v_n$ are on the cycles of $\overline{G}_{33}$ then we enter Case 1.1 or Case 1.2,
and otherwise repeating this process we can determine that the result graph $\widetilde{G}$ is $\theta_4$ or $\infty_5$ with trees attached.
If $\widetilde{G}\notin \mathcal{H}$,
then making some steps of $T_2$-transformation for $\widetilde{G}$ we can get the graph $\widetilde{G}_{T_2}$
which is isomorphic to one of the graphs $H_5$, $H_6$ and $H_7$ as in Fig. 2.

\textbf{Case 1.4.} Both $v_1$ and $v_n$ are on the same tree $T$.

Suppose that the tree $T$ is attached at the vertex $v_{l_m}$ of the cycle $C_k$.
Now we observe the path $v_n v_1v_{l_1}v_{l_2}\cdots v_{l_m}$.
Delete $v_{l_1}v_{l_2}$, and add the edge $v_1v_{l_2}$ if $(X_1+ X_{l_2} )^2\geq (X_n+ X_{l_2} )^2$,
and otherwise add the edge $v_nv_{l_2}$.
We make some steps of the above transformation until $v_{l_1}$ is on the cycle $v_{l_1}v_{s_2}\cdots v_{s_k}v_{l_1}$.
Next delete $v_{l_1}v_{s_2}$, and add the edge $v_1v_{s_2}$ if $(X_1+ X_{s_2} )^2\geq (X_n+ X_{s_2} )^2$,
and otherwise add the edge $v_nv_{s_2}$.
Then we can get the result graph $\overline{G}_{41}$ such that
$v_1$ and $v_n$ are on the same cycle of $\overline{G}_{41}$ or one of $v_1$ and $v_n$ is on the cycle and another is on the tree.
If the former rises then we enter Case 1.1, and otherwise we enter Case 1.3.

\textbf{Case 2.} $d_{G_{T_1}}(v_1,v_n) = 2$.

\textbf{Case 2.1.} $v_1$ and $v_n$ are on the same cycle $C_{t+3}$.

If $t\geq 2$ then set $C_{t+3}=v_1v_tv_n v_{l_1}v_{l_2}\cdots v_{l_t}v_1$.
Delete $v_{l_1}v_{l_2}$, and add the edge $v_1v_{l_1}$ if $(X_1+ X_{l_1} )^2\geq (X_n+ X_{l_2} )^2$,
and otherwise add the edge $v_nv_{l_2}$.
Repeat this process until the result graph $\widehat{G}_{11}$ contains the cycle $C_4 =v_1v_t v_n v_{l_1} v_1$.
Set $v_{s_1}v_{s_2}\cdots v_{s_k}v_{s_1}$ to be another cycle which is contained in $\widehat{G}_{11}$.
Delete $v_{s_1}v_{s_2}$, and add the edge $v_1v_{s_2}$ if $(X_1+X_{s_2} )^2\geq (X_n+X_{s_1} )^2$, and otherwise add the edge $v_nv_{s_1}$.
Repeating this process we can determine that the result graph $\widetilde{G}$ is $\theta_5$ or $\infty_6$ with trees attached.
If $\widetilde{G}\notin \mathcal{H}$,
then we make some steps of $T_2$-transformation for $\widetilde{G}$ until the result graph $\widetilde{G}_{T_2}$
is isomorphic to one of the graphs $G_1(p,q)$, $G_2(p,q)$ and $G_3(p,q)$ as in Fig. 1.

\textbf{Case 2.2.} $v_1$ and $v_n$ are on the different cycle.

In this case $G_{T_1}$ is the $b$-graph or $\infty$-graph with trees attached.
Suppose $G_{T_1}$ contains the cycle $C_{t+1}=v_1v_{l_1}v_{l_2}\cdots v_{l_t}v_1$.
If $t\geq 3$ then delete $v_{l_1}v_{l_2}$, and add the edge $v_1v_{l_2}$ if $(X_1+X_{l_2} )^2\geq (X_n+X_{l_1} )^2$,
and otherwise add the edge $v_nv_{l_1}$.
We denote by $\widehat{G}_{21}$ the result graph.
If $v_1$ and $v_n$ are on the same cycle of $\widehat{G}_{21}$ then we enter Case 2.1,
and otherwise repeat this process until the result graph $\widehat{G}_{22}$ contains the cycle $C_3$.
Let $v_n v_{s_1}v_{s_2}\cdots v_{s_k}v_n$ be another cycle of $\widehat{G}_{22}$.
Delete $v_{s_1}v_{s_2}$, and add the edge $v_1v_{s_1}$ if $(X_1+X_{s_1} )^2\geq (X_n+X_{s_2} )^2$, and otherwise add the edge $v_nv_{s_2}$.
We denote by $\widehat{G}_{23}$ the result graph.
If $v_1$ and $v_n$ are on the same cycle of $\widehat{G}_{23}$ then we enter Case 2.1,
and otherwise repeat this process until result graph
$\widetilde{G}$ is $\infty_5$, $b_7$ or $b_6$  with trees attached.
If $\widetilde{G}\notin \mathcal{H}$,
then we make some steps of $T_2$-transformation for $\widetilde{G}$
until the result graph $\widetilde{G}_{T_2}$ is isomorphic to one of the graphs $G_4(p,q)$, $G_5(p,q)$ and $G_6(p,q)$ as in Fig. 1.

\textbf{Case 2.3.} One of $v_1$ and $v_n$ is on the cycle and another is on the tree.

Suppose without loss of generality that $v_1$ is on the cycle and that $v_n$ is on the tree.
Then $G_{T_1}$ contains the cycle $C_{t+1}=v_1v_{l_1}v_{l_2}\cdots v_{l_t}v_1$.
If $t\geq 3$ then delete $v_{l_1}v_{l_2}$, and add the edge $v_1v_{l_2}$ if $(X_1+X_{l_2} )^2\geq (X_n+X_{l_1} )^2$,
and otherwise add the edge $v_nv_{l_1}$.
We denote by $\widehat{G}_{31}$ the result graph.
If $v_1$ and $v_n$ are on the same cycle of $\widehat{G}_{31}$ then we enter Case 2.1,
and otherwise repeat this process until the result graph $\widehat{G}_{32}$ contains the cycle $C_3=v_1v_{l_1}v_{l_2}v_1$.
Set $v_{s_1}v_{s_2}\cdots v_{s_k}v_{s_1}$ to be another cycle of $\widehat{G}_{32}$.
Delete $v_{s_1}v_{s_2}$, and add the edge $v_1v_{s_2}$ if $(X_1+X_{s_2} )^2\geq (X_n+X_{s_1} )^2$, and otherwise add the edge $v_nv_{s_1}$.
We denote by $\widehat{G}_{33}$ the result graph.
If $v_1$ and $v_n$ are on the cycles of $\widehat{G}_{33}$ then we enter Case 2.1 or Case 2.2,
and otherwise repeat this process until the result graph $\widetilde{G}$ is $\theta_4$ or $\infty_5$ with trees attached.
If $\widetilde{G}\notin \mathcal{H}$,
then we make some steps of $T_2$-transformation for $\widetilde{G}$ until the result graph $\widetilde{G}_{T_2}$
is isomorphic to one of graphs $G_7(p,q)$, $G_8(p,q)$, $G_9(p,q)$, $G_{10}(p,q)$, $G_{11}(p,q)$ and $G_{12}(p,q)$ as in Fig. 1.

\textbf{Case 2.4.} Both $v_1$ and $v_n$ are on the same tree $T$.

Suppose that the tree $T$ is attached at the vertex $v_{l_m}$ of the cycle $C_k$.
Let $v_t$ be the common neighbour of $v_1$ and $v_n$.
Now we distinguish two cases to discuss.

\textbf{Case 2.4.1.} Suppose that there exists the path $v_nv_tv_1v_{l_1}v_{l_2}\cdots v_{l_m}$.

Delete $v_{l_1}v_{l_2}$, and add the edge $v_1v_{l_2}$ if $(X_1+ X_{l_2} )^2\geq (X_n+ X_{l_2} )^2$, and otherwise add the edge $v_nv_{l_2}$.
We make some steps of the above transformation
until we obtain the cycle $v_{l_1}v_{s_2}v_{s_3}\cdots v_{s_k}v_{l_1}$.
Delete $v_{l_1}v_{s_2}$, and add the edge $v_1v_{s_2}$ if $(X_1+ X_{s_2} )^2\geq (X_n+ X_{s_2} )^2$,
and otherwise add the edge $v_nv_{s_2}$.
We denote by $\widehat{G}_{41}$ the result graph.
It can be observed that $v_1$ is on the cycle and $v_n$ is on the tree of $\widehat{G}_{41}$ or $v_1$ and $v_n$ are on the same cycle of $\widehat{G}_{41}$.
If the former rises then we enter Case 2.1, and otherwise we enter Case 2.3.

\textbf{Case 2.4.2.} Suppose that there exists the path $v_1v_tv_{l_1}v_{l_2}\cdots v_{l_m}$.

Delete $v_tv_{l_1}$, and add the edge $v_1v_{l_1}$
if $(X_1+ X_{l_1} )^2\geq (X_n+ X_{l_1} )^2$,
and otherwise add the edge $v_nv_{l_1}$.
We denote by $\widehat{G}_{42}$ the result graph.
If we obtain the path $v_nv_tv_1v_{l_1}v_{l_2}\cdots v_{l_m}$,
then we enter Case 2.4.1,
and otherwise repeat this process until we obtain the cycle $v_tv_{s_2}v_{s_3}\cdots v_{s_k}v_t$.
Delete $v_tv_{l_1}$, and add the edge $v_1v_{l_1}$ if $(X_1+ X_{l_1} )^2\geq (X_n+ X_{l_1} )^2$,
and otherwise add the edge $v_nv_{l_1}$.
We denote by $\widehat{G}_{43}$ the result graph.
It can be observed that one of $v_1$ and $v_n$ is on the cycle and another is on the tree of $\widehat{G}_{43}$,
and then we enter Case 2.3.

From the above argument we observe that
for any $G^c\in \mathcal {C}_n^c$ ($n\geq 12$),
there is a graph $H\in \mathcal{H}^{'}\cup \mathcal{H}^{''}$
such that $$X^{T}Q(G)X=\sum\limits_{v_iv_j\in E(G)}(X_i+X_j )^2\leq
\sum\limits_{v_iv_j\in E(H)}(X_i+X_j )^2=X^{T}Q(H)X.  \eqno{(2.4)}$$

Note that $Q(G^c)=(n-2)I_n+J_n-Q(G)$,
where $J_n$ denotes the all ones square matrix of order $n$.
Therefore, using the equations (2.2) and (2.4) we obtain

\begin{equation}
\begin{split}
\lambda (G^c)&=X^{T}Q(G^c)X\\
&=X^{T}((n-2)I_n+J_n)X-X^{T}Q(G)X\\
&\geq X^{T}((n-2)I_n+J_n)X-X^{T}Q(H)X\\
&= X^{T}Q(H^c)X\geq \lambda (Q(H^c)).~~~~~~~~~~~~~  \Box
\nonumber
\end{split}
\end{equation}

{\bf Lemma 2.4.} {\it If $H\in \mathcal{H}^{''}$,
then there is a graph $H_*\in \mathcal{H}^{'}$
such that $\lambda(H^c)\geq \lambda({H_*}^c).$}

{\bf Proof.} Let $X=(X_1, X_2, \ldots, X_n)^T$ be the unit first signless Laplacian eigenvector of $H^c$
such that $X_1 \geq X_2 \geq\dots\geq X_n$. Then $X_1> 0$ and $X_n< 0$.
We suppose without loss of generality that
$\vert X_n\vert \geq \vert X_1\vert$
for otherwise we can replace $X$ with $-X$.

For the graphs $H_i\in \mathcal{H}^{''}$ ($i=1, 2$),
deleting the edge $v_1v_n$, and adding the edge $v_nv_t$,
we get the result graph $H_{*}$ which is isomorphic to $G_1(p,q)$ or $G_2(p,q)$ in $\mathcal{H}^{'}$.
Thus we obtain

\begin{equation}
\begin{split}
X^{T}Q(H_{*})X-X^{T}Q(H)X&=(X_n+X_t )^2-(X_n+X_1)^2\\
&=(2X_n+X_t+X_1)(X_t-X_1)\geq 0.~~~~~~~~~~~~~~~~~~~~~~(2.5)
\nonumber
\end{split}
\end{equation}

For the graph  $H_i\in \mathcal{H}^{''}$ $(3\leq i\leq 7)$,
deleting the edge $v_1v_n$, and adding the edge $v_nv_k$,
we get the result graph $H_{*}$ which is isomorphic to $G_{\ell}(p,q)$ ($\ell=2, 6, 7, 8$ or $11$).
Thus we obtain

\begin{equation}
\begin{split}
X^{T}Q(H_{*})X-X^{T}Q(H)X&=(X_n+X_k )^2-(X_1+X_n)^2\\
&=(2X_n+X_k+X_1)(X_k-X_1)\geq 0~~~(2.6).
\nonumber
\end{split}
\end{equation}

Therefore, combining the equations (2.2), (2.5) and (2.6) we have\vskip 2mm

\begin{equation}
\begin{split}
\lambda({H}^c)&=X^{T}Q({H}^c)X\\
&=X^{T}((n-2)I_n+J_n)X-X^{T}Q(H)X\\
&\geq X^{T}((n-2)I_n+J_n)X-X^{T}Q(H_{*})X\\
&= X^{T}Q(H_{*}^c)X\geq \lambda (H_{*}^c).~~~~~~~~~\Box
\nonumber
\end{split}
\end{equation}

{\bf Lemma 2.5} \cite{SCLi}. {\it Let $G$ be a simple graph.
Then $\lambda(G) \leq \delta (G)$, where $\delta(G)= min\{d_{G}(v),v\in V(G) \}$.} \vskip 3mm

{\bf Lemma 2.6.} {\it Let $p$ and $q$ be two positive integers such that $p+q=n-5$ ($n\geq 12$).
Then $\lambda(G^c_{1}(p,q)) > \lambda(G^c_{1}(n-5, 0)).$}

{\bf Proof.} We without loss of generality assume $p\geq q\geq 1$.
Next we will prove $\lambda(G^c_{1}(p,q)) > \lambda(G^c_{1}(p+1,q-1))$.
Suppose $X=(X_1, X_2, \ldots, X_n)^T$ is the unit first signless Laplacian eigenvector of $G^c_{1}(p,q)$,
and let $k_1$=$\lambda(G^c_{1}(p,q))$.
Then, by the symmetry of $G^c_{1}(p,q)$ and the signless Laplacian eigen-equation (2.3),
we have
\begin{flushleft}
$\left\{
\begin{array}{ll}
\ (k_1-(q+1))X_1 =X_6+q X_7,& \\
\ (k_1-(p+q+3))X_2 =(p-1)X_2+X_3+X_4+X_5+X_6+q X_7,& \\
\ (k_1-(p+q+2))X_3 =pX_2+X_4+X_5+q X_7,& \\
\ (k_1-(p+q+2))X_4=pX_2+X_3+X_5+q X_7,&   \\
\ (k_1-(p+q+2))X_5 =pX_2+X_3+X_4+q X_7,& \\
\ (k_1-(p+1))X_6=X_1+p X_2,& \\
\ (k_1-(p+q+3))X_7 =X_1+p X_2+X_3+X_4+X_5+(q-1)X_7.&  \hbox{}
\end{array}
\right.$
\end{flushleft}
We can transform the above equations into a matrix equation $(k_1I_7-Q_1)X^\prime =0$,
where $X^\prime =(X_1,X_2,X_3,X_4,X_5,X_6,X_7)^T$ and

$$Q_1=
\begin{pmatrix}
q+1 &0&0&0&0&1& q \\
0&2p+q+2& 1& 1& 1& 1& q \\
0&p& p+q+2& 1& 1& 0& q \\
0&p& 1& p+q+2& 1& 0&q \\
0&p& 1& 1& p+q+2& 0&q \\
1& p& 0&0&0&p+1&0& \\
1& p & 1& 1& 1& 0& p+2q+2 &
\end{pmatrix}
$$
Let $f_{1}(x;p,q)=det(xI_7-Q_1)$. Then we have

$f_{1}(x;p,q )=x^7-(7p+7q+12)x^6-(-20p^2-41pq-67p-20q^2-67q-57)x^5
-(30p^3+96p^2q+150p^2+96pq^2+306pq+248p+30q^3+150q^2+248q+138)x^4
-(-25p^4-114p^3q-170p^3-178p^2q^2-540p^2q-416p^2-114pq^3-540pq^2
-844pq-447p-25q^4-170q^3-416q^2-447q-180)x^3-(11p^5+71p^4q+100p^4
+158p^3q^2+454p^3q+330p^3+158p^2q^3+708p^2q^2+1038p^2q+518p^2+71pq^4
+454pq^3+1038pq^2+1046pq+397p+11q^5+100q^4+330q^3+518q^2
+397q+120)x^2-(-2p^6-21p^5q-27p^5-66p^4q^2-177p^4q-119p^4-94p^3q^3
-396p^3q^2-536p^3q-247p^3-66p^2q^4-396p^2q^3-834p^2q^2-771p^2q
-267p^2-21pq^5-177pq^4-536pq^3-771pq^2-537pq-146p-2q^6-27q^5
-119q^4-247q^3-267q^2-146q-32)x-2p^6q-2p^6-10p^5q^2-24p^5q-14p^5
-20p^4q^3-78p^4q^2-94p^4q-38p^4-20p^3q^4-112p^3q^3-212p^3q^2
-172p^3q-50p^3-10p^2q^5-78p^2q^4-212p^2q^3-268p^2q^2-156p^2q-32p^2
-2pq^6-24pq^5-94pq^4-172pq^3-156pq^2-64pq-8p-2q^6-14q^5-38q^4-50q^3-32q^2-8q.$

Therefore, we obtain
\begin{equation}
\begin{split}
&f_{1}(x;p,q )-f_{1}(x;p+1,q-1)\\
&=(p-q+1)(-q-3-p+x)(-2q-2p+x)(x-p-q-1)^3~~~~~~~~~~~~~~(2.7)
\nonumber
\end{split}
\end{equation}

Note that $G^c_{1}(p,q)$ is connected and $\delta(G^c_{1}(p,q))=q+1$.
By Lemma 2.5, we have $0<k_1\leq q+1$.
Since $f_{1}(k_1;p,q)=0$ and $p\geq q\geq 1$,
it follows that $f_{1}(k_1;p+1,q-1)>0$ from the equation (2.7).
This shows that $\lambda(G^c_{1}(p,q))>\lambda(G^c_{1}(p+1,q-1))$,
which implies $\lambda(G^c_{1}(p,q)) > \lambda(G^c_{1}(n-5,0)).$    $\Box$ \vskip 3mm

{\bf Lemma 2.7.} {\it Let $p$ and $q$ be two positive integers such that $p+q=n-5$ ($n\geq 12$).
Then $\lambda(G^c_{2}(p,q)) > \lambda(G^c_{2}(n-5,0))$.}

{\bf Proof.} Suppose $X=(X_1, X_2, \cdots, X_n)^T$ is the unit first signless Laplacian eigenvector of $ G^c_{2}(p,q) $,
and let $k_2$=$\lambda(G^c_{2}(p,q))$.
Then, by the symmetry of $G^c_{2}(p,q)$ and the signless Laplacian eigen-equation (2.3),
we have

\begin{flushleft}
$\left\{
\begin{array}{ll}
\ (k_2-(q+1))X_1 =X_6+q X_7,& \\
\ (k_2-(p+q+3))X_2 =(p-1)X_2+X_3+X_4+X_5+X_6+q X_7,& \\
\ (k_2-(p+q+2))X_3 =pX_2+X_4+X_5+q X_7,& \\
\ (k_2-(p+q+2))X_4=pX_2+X_3+X_6+q X_7,&  \\
\ (k_2-(p+q+1))X_5 =pX_2+X_3+q X_7,& \\
\ (k_2-(p+2))X_6=X_1+p X_2+X_4,& \\
\ (k_2-(p+q+3))X_7 =X_1+pX_2+X_3+X_4+X_5+(q-1)X_7.&  \hbox{}
\end{array}
\right.$
\end{flushleft}
We can transform the above equations into a matrix equation $(k_2I_7-Q_2)X^\prime =0$,
where $X^\prime =(X_1,X_2,X_3,X_4,X_5,X_6,X_7)^T$ and

$$Q_2=
\begin{pmatrix}
q+1 &0&0&0&0&1& q \\
0&2p+q+2& 1& 1& 1& 1& q \\
0&p& p+q+2& 1& 1& 0& q \\
0&p& 1& p+q+2& 0& 1&q \\
0&p& 1& 0& p+q+1& 0&q \\
1& p& 0&1&0&p+2&0& \\
1& p & 1& 1& 1& 0& p+2q+2 &
\end{pmatrix}
$$
Let $f_{2}(x;p,q)=det(xI_7-Q_2)$. Then we have

$f_{2}(x;p,q)=x^7-(7p+7q+12)x^6-(-20p^2-41pq-67p-20q^2-68q-57)x^5-(30p^3+96p^2q+150p^2+96pq^2+311pq+248p+30q^3+156q^2+253q+136)x^4
-(-25p^4-114p^3q-170p^3-178p^2q^2-549p^2q-416p^2-114pq^3-563pq^2-867pq-444p-25q^4-184q^3-441q^2-449q-169)x^3
-(11p^5+71p^4q+100p^4+158p^3q^2+461p^3q+330p^3+158p^2q^3+738p^2q^2+1071p^2q+520p^2+71pq^4+493pq^3+1116pq^2+1079pq+389p
+11q^5+116q^4+375q^3+544q^2+380q+100)x^2-(-2p^6-21p^5q-27p^5-66p^4q^2-179p^4q-119p^4-94p^3q^3-411p^3q^2-553p^3q-252p^3
-66p^2q^4-429p^2q^3-903p^2q^2-818p^2q-280p^2-21pq^5-206pq^4-623pq^3-851pq^2-566pq-146p-2q^6-36q^5-154q^4-285q^3-269q^2-128q-20)x
-2p^6q-2p^6-10p^5q^2-24p^5q-14p^5-20p^4q^3-80p^4q^2-96p^4q-40p^4-20p^3q^4-120p^3q^3-228p^3q^2-188p^3q-60p^3-10p^2q^5-90p^2q^4
-248p^2q^3-310p^2q^2-192p^2q-48p^2-2pq^6-32pq^5-126pq^4-216pq^3-190pq^2-88pq-16p-4q^6-24q^5-54q^4-58q^3-32q^2-8q.$

Therefore, we obtain
\begin{equation}
\begin{split}
&f_{2}(x;n-5-q,q)-f_{2}(x;n-5,0)=-q(n-4-x)g_1(x)~~~~~~~~~~~~~~~~~~~~~~~~~~~~~(2.8)
\nonumber
\end{split}
\end{equation}
where $g_1(x)=(n-q-4)x^4+(-5n^2+5nq+41n-21q-84)x^3+(9n^3-9n^2q-110n^2+74nq+447n-151q-609)x^2
+(-7n^4+7n^3q+113n^3-85n^2q-679n^2+339nq+1799n-436q-1760)x+2n^5-2n^4q-40n^4+32n^3q+316n^3-188n^2q-1228n^2+474nq+2330n-428q-1708.$

Note that $G^c_{2}(p,q)$ is connected and $\delta(G^c_{2}(p,q))=min\{q+1, p+2\}$.
By Lemma 2.5, we have $0<k_2\leq min\{q+1, p+2\}.$
Claim A in the Appendix shows that $g_1(x)>0$ when $0<x\leq min\{q+1, p+2\}.$
Since $f_{2}(k_2;n-5-q,q)=0$,
it follows that $f_{2}(k_2;n-5,0)>0$ from the equation (2.8).
This shows that $\lambda(G^c_{2}(p,q))>\lambda(G^c_{2}(n-5,0))$.  ~~~~~$\Box$\vskip 3mm

For $G_{k}(p,q)\in \mathcal{H}^{'}$ $(k=3, 6, 7, 8, 9, 10, 11, 12)$,
we can prove as the proof of Lemma 2.7 that the following Lemma is true.\vskip 2mm

{\bf Lemma 2.8.} {\it  Let $p$ and $q$ be two positive integers and $n\geq 12$.
Then we have
\begin{enumerate}
\setlength{\parskip}{0ex}
\item[\rm (i)]
when $p+q=n-5$, $\lambda(G^c_{s}(p,q))>\lambda(G^c_{s}(n-5,0))~~(s=8, 9, 11)$,
\item[\rm (ii)]
when $p+q=n-6$, $\lambda(G^c_{t}(p,q))>\lambda(G^c_{t}(n-6,0))~~(t=3, 7, 10)$, and\\
  $\lambda(G^c_{6}(p,q))>\lambda(G^c_{6}(0,n-6))$,
\item[\rm (iii)]
when $p+q=n-7$, $\lambda(G^c_{12}(p,q))>\lambda(G^c_{12}(n-7,0))$.
\end{enumerate} }

{\bf Lemma 2.9.} {\it Let $p$ and $q$ be two positive integers such that $p+q=n-7$ ($n\geq 12$).
Then $\lambda(G^c_{4}(p,q)) > \lambda(G^c_{4}(n-7,0)).$ }

{\bf Proof.}  Without loss of generality assume $p\geq q\geq 1$.
We will first prove $\lambda(G^c_{4}(p,q)) > \lambda(G^c_{4}(p+1,q-1))$.
Suppose $X=(X_1, X_2, \ldots, X_n)^T$ is the unit first signless Laplacian eigenvector of $ G^c_{4}(p,q) $,
and let $k_4$=$\lambda(G^c_{4}(p,q))$.
Then, by the symmetry of $G^c_{4}(p,q)$ and the signless Laplacian eigen-equation (2.3),
we have

\begin{flushleft}
$\left\{
\begin{array}{ll}
\ (k_4-(q+3))X_1 =2X_5+X_6+q X_7,& \\
\ (k_4-(p+q+5))X_2 =(p-1)X_2+2X_3+X_4+2X_5+X_6+q X_7,& \\
\ (k_4-(p+q+4))X_3 =pX_2+X_4+2X_5+X_6+q X_7,& \\
\ (k_4-(p+q+4))X_4=pX_2+2X_3+2X_5+q X_7,&  \\
\ (k_4-(p+q+4))X_5 =X_1+pX_2+2X_3+X_4+q X_7,& \\
\ (k_4-(p+3))X_6=X_1+p X_2+2X_3,& \\
\ (k_4-(p+q+5))X_7 =X_1+p X_2+2X_3+X_4+2X_5+(q-1)X_7.&  \hbox{}
\end{array}
\right.$
\end{flushleft}
We can transform the above equations into a matrix equation $(k_4I_7-Q_4)X^\prime =0$,
where $X^\prime =(X_1,X_2,X_3,X_4,X_5,X_6,X_7)^T$ and

$$Q_4=
\begin{pmatrix}
q+3 &0&0&0&2&1& q \\
0&2p+q+4& 2& 1& 2& 1& q \\
0&p& p+q+4& 1& 2& 1& q \\
0&p& 2& p+q+4& 2& 0&q \\
1&p& 2& 1& p+q+4& 0&q \\
1& p& 2&0&0&p+3&0& \\
1& p & 2& 1& 2& 0& p+2q+4 &
\end{pmatrix}
$$
Let $f_{4}(x;p,q)=det(xI_7-Q_4)$. Then we have

$f_{4}(x;p,q)=x^7-(7p+7q+26)x^6-(-20p^2-41pq-149p-20q^2-149q-276)x^5-(30p^3+96p^2q+342p^2+96pq^2+700pq+1269p+30q^3+342q^2+1269q+1548)x^4
-(-25p^4-114p^3q-398p^3-178p^2q^2-1270p^2q-2253p^2-114pq^3-1270pq^2-4613pq-5506p-25q^4-398q^3-2253q^2-5506q-4916)x^3
-(11p^5+71p^4q+242p^4+158p^3q^2+1100p^3q+1905p^3+158p^2q^3+1716p^2q^2+6087p^2q+7092p^2+71pq^4+1100pq^3+6087pq^2
+14552pq+12696p+11q^5+242q^4+1905q^3+7092q^2+12696q+8688)x^2
-(-2p^6-21p^5q-69p^5-66p^4q^2-445p^4q-747p^4-94p^3q^3-990p^3q^2-3411p^3q-3858p^3-66p^2q^4-990p^2q^3-5328p^2q^2
-12402p^2q-10516p^2-21pq^5-445pq^4-3411pq^3-12402pq^2-21668pq-14432p-2q^6-69q^5-747q^4-3858q^3-10516q^2-14432q-7616)x
-2p^6q-6p^6-10p^5q^2-64p^5q-102p^5-20p^4q^3-202p^4q^2-668p^4q-724p^4-20p^3q^4-288p^3q^3-1494p^3q^2-3356p^3q
-2736p^3-10p^2q^5-202p^2q^4-1494p^2q^3-5264p^2q^2-8904p^2q-5712p^2-2pq^6-64pq^5-668pq^4-3356pq^3-8904pq^2-11856pq
-6016p-6q^6-102q^5-724q^4-2736q^3-5712q^2-6016q-2304.$

Therefore, we obtain
\begin{equation}
\begin{split}
&f_{4}(x;p,q )-f_{4}(x;p+1,q-1)\\
&=-(p-q+1)(-x+3+q+p)(-x+2+q+p)g_2(x)~~~~~~~~~~~~~~~~~~~~~~~(2.9)
\nonumber
\end{split}
\end{equation}
where $g_2(x)=-x^3+(4p+4q+11)x^2+(-5p^2-10pq-5q^2-29p-29q-46)x+2p^3+6p^2q+6pq^2+2q^3+18p^2+36pq+18q^2+56p+56q+72.$

Note that $G^c_{4}(p,q)$ is connected and $\delta(G^c_{4}(p,q))=q+3$.
By Lemma 2.5, we have $0<k_4\leq q+3$.
Claim B in the Appendix shows that $g_2(x)>0$ when $0<x\leq q+3.$
Since $f_{4}(k_4;p,q)=0$,
it follows $f_{4}(k_4;p+1,q-1)>0$ from the equation (2.9).
This shows that $\lambda(G^c_{4}(p,q))>\lambda(G^c_{4}(p+1,q-1))$,
from which we obtain $\lambda(G^c_{4}(p,q)) > \lambda(G^c_{4}(n-7,0)).$ $\Box$ \vskip 3mm

For $G_{5}(p,q)\in \mathcal{H}^{'}$,
we can prove as the proof of Lemma 2.9 that the following result is true.\vskip 2mm

{\bf Lemma 2.10.} {\it  Let $p$ and $q$ be two positive integers such that $p+q=n-5$ $(n\geq 12)$.
Then $\lambda(G^c_{5}(p,q)) > \lambda(G^c_{5}(n-5,0)).$}\vskip 3mm

{\bf Theorem 2.11.} {\it For any graph  $G^c\in \mathcal {C}_n^c$ $(n\geq 12)$,
we have $\lambda(G^c)\geq \lambda(G^c_{1}(n-5,0))$
with equality if and only if $G\cong G_{1}(n-5,0)$.}

{\bf Proof.} Let $f_{1}(x;p,q)$, $f_{2}(x;p,q)$ and $f_{4}(x;p,q)$
be as in the proof of Lemmas 2.6, 2.7 and 2.9, respectively.

Let $g_{12}(x)=f_{1}(x;n-5,0)-f_{2}(x;n-5,0)$.
Then $g_{12}(x)=(n-x-3)^2(-2n^2+(-x-18)n-2x^2+8x+40)$.
Since $0<x\leq 1$ and $n\geq 12$,
we have $g_{12}(x)>0$.
Lemma 2.5 shows $0<\lambda(G^c_{1}(n-5,0)), \lambda(G^c_{2}(n-5,0))\leq 1$.
Note that $f_{1}(\lambda(G^c_{1}(n-5,0));n-5,0)=0$.
Thus, we obtain $f_{2}(\lambda(G^c_{1}(n-5,0));n-7,0)<0$.
This implies $\lambda(G^c_{1}(n-5,0))<\lambda(G^c_{2}(n-7,0)).$

Let $g_{14}(x)=f_{1}(x;n-5,0)-f_{4}(x;n-7,0).$
Then $g_{14}(x)=(n-3-x)g_4(x)$,
where $g_4(x)=(2n-9)x^4+(-10n^2+92n-208)x^3+(18n^3-249n^2+1149n-1814)x^2+(-14n^4+258n^3-1791n^2+5607n-6672)x
+4n^5-92n^4+850n^3-3956n^2+9242n-8592.$

Claim C ~in the Appendix shows that $g_{4}(x)>0$, and so $g_{14}(x)>0$ when $0< x\leq 1$ and $n\geq 12$.
Lemma 2.5 shows $0<\lambda(G^c_{1}(n-5,0)), \lambda(G^c_{4}(n-7,0))\leq 1$.
Note that $f_{1}(\lambda(G^c_{1}(n-5,0));n-5,0)=0$.
Thus, we obtain $f_{4}(\lambda(G^c_{1}(n-5,0));n-7,0)<0$.
This implies $\lambda(G^c_{1}(n-5,0))<\lambda(G^c_{4}(n-7,0)).$

Similarly, we can verify that the following results are true.

\begin{equation}
\begin{split}
&\lambda(G^c_{1}(n-5,0))<\lambda(G^c_{s}(n-5,0))~~(s=5, 8, 9, 11),\\
&\lambda(G^c_{1}(n-5,0))<\lambda(G^c_{t}(n-6,0))~~(t=3, 7, 10),\\
&\lambda(G^c_{1}(n-5,0))<\lambda(G^c_{6}(0,n-6))\\
&\lambda(G^c_{1}(n-5,0))<\lambda(G^c_{12}(n-7,0))
\nonumber
\end{split}
\end{equation}

Combining these inequations with Lemmas 2.3, 2.4, 2.6, 2.7, 2.8, 2.9 and 2.10,
we determine that the result is true.  $\Box$

\appendix
\appendixpage
\begin{appendices}

~~~~{\bf Claim A.}  {\it Let $g_1(x)$ be as in the proof of Lemma 2.7.
Then $g_1(x)>0$ when $0<x\leq min\{q+1, p+2\}$.}

{\bf Proof.}  We take the derivatives for $g_1(x)$ as follows.

${g_1}'(x)=4(n-q-4)x^3+3(-5n^2+5nq+41n-21q-84)x^2+2(9n^3-9n^2q-110n^2+74nq+447n-151q-609)x
-7n^4+7n^3q+113n^3-85n^2q-679n^2+339nq+1799n-436q-1760.$

${g_1}''(x)=12(n-q-4)x^2+6(-5n^2+5nq+41n-21q-84)x+18n^3-18n^2q-220n^2+148nq+894n-302q-1218.$

${g_1}'''(x)=24(n-q-4)x-30n^2+30nq+246n-126q-504.$

Then ${g_1}'''(x)$ is monotonously increasing for $x$.

Suppose $min\{q+1,p+2\}=p+2$.
Then $0<x\leq p+2$.
Since $p\geq 0$ and $n=p+q+5 \geq 12$, we have $n\geq 2p+6$ and $q\geq4$,
and so ${g_1}'''(x)\leq {g_1}'''(p+2)=-6(n-q-4)(5n-4p-29)<0$.
This shows ${g_1}''(x)$ is monotonously decreasing for $x$, and so
${g_1}''(x)\geq {g_1}''(p+2)=18n^3+(-30p-18q-280)n^2+(12p^2+30pq+294p+208q+1434)n-12p^2q-48p^2-174pq-696p-602q-2418>0$.

Thus, ${g_1}'(x)$ is monotonously increasing for $x$, and so
${g_1}'(x) \leq {g_1}'(p+2)=-7n^4+(18p+7q+149)n^3+(-15p^2-18pq-280p-121q-1179)n^2+(4p^3+15p^2q+147p^2+208pq+1434p+695q+4111)n
-4p^3q-16p^3-87p^2q-348p^2-602pq-2418p-1324q-5332<0$.

We determine that $g_1(x)$ is monotonously decreasing for $x$, and so

$g_1(x)\geq g_1(p+2)=2n^5+(-7p-2q-54)n^4+(9p^2+7pq+149p+46q+578)n^3+(-5p^3-9p^2q-140p^2-121pq-1179p-394q-3066)n^2
+(p^4+5p^3q+49p^3+104p^2q+717p^2+695pq+4111p+1488q+8060)n-p^4q-4p^4-29p^3q-116p^3-301p^2q-1209p^2-1324pq-5332p-2088q-8400>0.$

If $min\{q+1, p+2\}=q+1$
then $0<x\leq q+1$.
We can similarly prove $g_1(x)>0$.    $\Box$ \vskip 2mm

{\bf Claim B.}  {\it Let $g_2(x)$ be as in the proof of Lemma 2.9.
Then $g_2(x)>0$ when $0<x\leq q+3$.}

{\bf Proof.}  We take the derivatives for $g_2(x)$ as follows.

${g_2}'(x)=-3x^2+2(4p+4q+11)x-5p^2-10pq-5q^2-29p-29q-46.$

${g_2}''(x)=-6x+8p+8q+22.$

We observe that ${g_2}''(x)$ is monotonously decreasing for $x$,
and so ${g_2}''(x)\geq {g_2}''(q+3)=2q+4+8p>0.$
This shows ${g_2}'(x)$ is monotonously increasing for $x$, and so
${g_2}'(x) \leq {g_2}'(q+3)=(-2p-1)q-5p^2-5p-7<0.$
We determine that $g_2(x)$ is monotonously decreasing for $x$, and so
$g_2(x) \geq g_2(q+3)=(p^2+p-2)q+2p^3+3p^2+5p+6>0.$     $\Box$ \vskip 3mm

{\bf Claim C.}  {\it Let $g_4(x)$ be as in the proof of Theorem 2.11.
Then $g_4(x)>0$ when $0<x\leq 1$.}

{\bf Proof.} Taking the derivatives for $g_4(x)$, we obtain

${g_4}'(x)=4(2n-9)x^3+3(-10n^2+92n-208)x^2+(2(18n^3-249n^2+1149n-1814))x-14n^4+258n^3-1791n^2+5607n-6672,$

${g_4}''(x)=12(2n-9)x^2+6(-10n^2+92n-208)x+36n^3-498n^2+2298n-3628,$

${g_4}'''(x)=24(2n-9)x-60n^2+552n-1248.$

Note that $0<x\leq 1$ and $n\geq 12$.
We have ${g_4}'''(x)$ is monotonously increasing for $x$,
and so ${g_4}'''(x)\leq {g_4}'''(1)=-60n^2+600n-1464<0.$

This shows ${g_4}''(x)$ is monotonously decreasing for $x$, and so
${g_4}''(x)\geq {g_4}''(1)=36n^3-558n^2+2874n-4984>0.$

Thus, ${g_4}'(x)$ is monotonously increasing for $x$, and so
${g_4}'(x)\leq {g_4}'(1)=-14n^4+294n^3-2319n^2+8189n-10960<0.$

We determine that $g_4(x)$ is monotonously decreasing for $x$, and so
$g_4(x)\geq g_4(1)=4n^5-106n^4+1126n^3-6006n^2-17295>0$.   $\Box$

\end{appendices}
\end{document}